\documentclass{amsart}

\usepackage{geometry}
\newtheorem{theorem}{Theorem}[section]
\newtheorem{example}[theorem]{Example}

\newcommand{\rank}{\text{rank}}

\begin{document}

\title{On a remarkable class of paracontact metric manifolds}

\author{Ver\'onica Mart\'{\i}n-Molina
 }

\address{Departamento de Matem\'aticas - IUMA, Universidad de Zaragoza, 
Pedro Cerbuna 12, 50009 Zaragoza, Spain}
\email{vmartin@unizar.es}

\begin{abstract}
We study a remarkable class of paracontact metric manifolds which have no contact metric counterpart: the paracontact metric $(-1,\widetilde\mu)$-spaces which are not paraSasakian (i.e. have $\widetilde h\neq0$). We present explicit examples with $\widetilde h$ of every possible constant rank and some with non-constant rank, which were not known to exist until recently.
\end{abstract}

\keywords{paracontact metric manifold; contact metric manifold; $(\kappa,\mu)$-spaces; nullity distribution; paraSasakian manifold.}

\maketitle

\section{Introduction}

Paracontact metric manifolds $(M,\widetilde\varphi,\xi,\eta,\widetilde g)$ have been studied by many authors in the recent years, particularly since the appearance of \cite{zamkovoy}. A special class among them is that of the $(\widetilde\kappa,\widetilde\mu)$-spaces, which satisfy \cite{CKM}
\begin{equation}\label{kappamu}
R(X,Y)\xi=\widetilde{\kappa}(\eta(Y)X-\eta(X)Y)+\widetilde{\mu}(\eta(Y) \widetilde{h} X-\eta(X) \widetilde{h} Y),
\end{equation}
for all $X,Y$ vector fields on $M$, where $\widetilde\kappa$ and $\widetilde\mu$ are constants and $\widetilde h=\frac12 L_\xi \widetilde\varphi$. These spaces include the paraSasakian manifolds \cite{kaneyuki,zamkovoy} and certain $g$-natural paracontact metric structures constructed on unit tangent sphere bundles \cite{mio-calvaruso}, among others.

Although the nullity condition~\eqref{kappamu} seems very technical, paracontact metric $(\widetilde\kappa,\widetilde\mu)$-spaces appear naturally when studying the relation between contact and paracontact metric geometry. Indeed, any non-Sasakian contact metric $(\kappa,\mu)$-space admits two paracontact metric $(\widetilde\kappa,\widetilde\mu)$-structures with the same contact form  and, under some natural conditions, every non-paraSasakian paracontact metric $(\widetilde{\kappa},\widetilde{\mu})$-space accepts a contact metric $(\kappa,\mu)$-structure with the same contact form \cite{mino-pacific,CKM}.

However, there are also some important differences between a contact metric $(\kappa,\mu)$-space $(M,\varphi,\xi,\eta,g)$ and a paracontact metric $(\widetilde\kappa,\widetilde \mu)$-space  $(M,\widetilde\varphi, \xi, \eta, \widetilde{g})$. First of all, while they satisfy  $h^2=(\kappa-1)\varphi^2$ and $\widetilde{h}^2=(\widetilde\kappa+1)\widetilde\varphi^2$, respectively, the first condition implies that  $\kappa \leq 1$ but the second one does not give any type of restriction over $\widetilde\kappa$ because $\widetilde{g}$ is not positive definite \cite{blair95,CKM}.

Another difference is that, in the contact metric case, $\kappa=1$ is also  equivalent to the manifold being Sasakian, i.e.  $h^2=0$ implies $h=0$. Nevertheless, there are non-paraSasakian paracontact metric $(-1,\widetilde\mu)$-spaces, i.e. with $\widetilde{h}^2=0$ but $\widetilde h\neq0$.

The first examples of these remarkable paracontact metric manifolds shown in the literature all had $\rank(\widetilde h)=n$ and $\widetilde{\mu}=0$ or $2$, \cite{CP,mino-pacific,CKM,murathan}. Indeed, until very recently, there seemed to be no literature discussing the rank of $\widetilde h$, if it had to be constant or why the values of $\mu$ zero and two were important.

This motivated the paper \cite{mio}, where the author presented a local classification of paracontact metric $(-1,\widetilde\mu)$-spaces in terms of the rank of $\widetilde h$, examples of paracontact metric $(-1,2)$-spaces with every possible constant rank of $\widetilde h$ and an explanation of why the values of $\mu$ zero and two are special. Later, the author also  wrote \cite{mio2}, where she gave an alternative proof of her main result, examples of paracontact metric $(-1,0)$-spaces with every possible constant rank of $\widetilde h$ and examples of paracontact metric $(-1,\widetilde\mu)$-spaces where $\widetilde h$ is of non-constant rank.

In the present paper, after the preliminaries section, we will summarize what is known about the remarkable class of paracontact metric $(-1,\widetilde\mu)$-spaces with $\widetilde h\neq0$, which have no contact metric counterpart.

\section{Preliminaries}

\emph{Almost paracontact manifolds} are $(2n+1)$-dimensional smooth manifolds endowed with a $(1,1)$-tensor $\widetilde\varphi$, a vector field $\xi$ and a $1$-form $\eta$ such that $\widetilde\varphi^2=I-\eta \otimes \xi$, $\eta(\xi)=1$ and $\widetilde\varphi$ induces a paracomplex structure on $\mathcal{D}=\text{ker} \, \eta$, i.e. the eigendistributions $\mathcal{D}^{\pm}$ corresponding to the eigenvalues $\pm 1$ of $\widetilde\varphi$ are both of dimension $n$ \cite{kaneyuki,zamkovoy}.

If the almost paracontact manifold admits a pseudo-Riemannian metric $\widetilde g$ of signature $(n+1,n)$ such that $\widetilde g(\widetilde\varphi X,\widetilde\varphi Y)=-\widetilde g(X,Y)+\eta(X)\eta(Y)$ and $d\eta(X,Y)=\widetilde g(X,\widetilde \varphi Y)$ for any vector fields $X$ and $Y$, then $M$ is called a \emph{paracontact metric manifold} and $(\widetilde\varphi,\xi,\eta,\widetilde g)$ its paracontact metric structure, \cite{zamkovoy}. We refer to \cite{mio-calvaruso2} for a recent survey on this type of manifold.

Given a a paracontact metric manifold, the tensor field $\widetilde{h}:=\frac{1}{2}L_{\xi} \widetilde\varphi$ is symmetric with respect to $\widetilde{g}$, i.e. $\widetilde g(\widetilde hX,Y)=\widetilde g(X,\widetilde hY)$, for all $X,Y$, and satisfies $\widetilde h \widetilde\varphi=-\widetilde\varphi \widetilde h$ and $\widetilde h\xi=\text{tr} \, \widetilde h=0$ \cite{zamkovoy}. Moreover, $\widetilde h=0$ if and only if $\xi$ is Killing, in which case the manifold is said to be a \emph{$K$-paracontact manifold}.

An almost paracontact structure is called \emph{normal} if and only if the tensor $[\widetilde{\varphi},\widetilde{\varphi}]-2d\eta\otimes\xi=0$, where $[\widetilde{\varphi},\widetilde{\varphi}]$ is the Nijenhuis tensor of $\widetilde{\varphi}$ \cite{zamkovoy}. A normal paracontact metric manifold is said to be a \emph{paraSasakian manifold} and is in particular $K$-paracontact. The converse holds in dimension $3$ \cite{calvaruso} and always for $(-1,\widetilde\mu)$-spaces \cite[Th.~3.1]{mio}. However, it is not true in general, \cite[Ex.~2.1]{mio}.

Every paraSasakian manifold satisfies
\begin{equation}\label{parasasakian-r}
R(X,Y)\xi =-(\eta(Y)X-\eta(X)Y),
\end{equation}
for every $X,Y$ on $M$. The converse is not true, since Examples~3.8--3.11 of \cite{mio} and Examples~4.1 and~4.5 of \cite{mio2} show that there are examples of paracontact metric manifolds satisfying Eq.~\eqref{parasasakian-r} but with $\widetilde{h}\neq0$ (and therefore not K-paracontact or paraSasakian).

\section{Classification and examples}

Many examples of paraSasakian manifolds are known. For instance, hyperboloids $\mathbb{H}^{2n+1}_{n+1}(1)$ and the hyperbolic Heisenberg group ${\mathcal H}^{2n+1}=\mathbb{R}^{2n}\times\mathbb{R}$, \cite{ivanov}. We can also obtain ($\eta$-Einstein) paraSasakian manifolds from contact metric $(\kappa,\mu)$-spaces with $|1-\frac{\mu}{2}|<\sqrt{1-\kappa}$. In particular, the tangent sphere bundle $T_1N$ of any space form $N(c)$ with $c<0$ admits a canonical $\eta$-Einstein paraSasakian structure,  \cite{nuestro-mino}. Finally, we can see how to construct explicitly a paraSasakian structure on a Lie group, \cite[Example 3.4]{mio}, or directly on the unit tangent sphere bundle, \cite{mio-calvaruso}.

On the other hand, until \cite{mio} and \cite{mio2} only the following examples of paracontact metric $(-1,\mu)$-spaces $(M^{2n+1},\widetilde{\varphi},\xi,\eta, \widetilde{g})$ with $\widetilde{h}\neq0$ were known (cited here in chronological order):
\begin{itemize}

\item paracontact metric $(-1,2)$-space with $\text{rank}(\widetilde h)=n=2$, \cite{mino-pacific}.

\item paracontact metric $(-1,2)$-spaces with $\text{rank}(\widetilde h)=n=\text{arbitrary}$,  \cite{CKM}.

\item paracontact metric $(-1,2)$-space with $\text{rank}(\widetilde h)=n=1$, \cite{murathan}.

\item paracontact metric $(-1,0)$-space with $\text{rank}(\widetilde h)=n=1$, \cite{CP}.

\end{itemize}


We will first show why there only seem to be examples of paracontact metric $(-1,\widetilde\mu)$-spaces with $\widetilde\mu=0$ or $\widetilde\mu=2$. Given a paracontact metric structure $(\widetilde{\varphi},\xi,\eta,\widetilde g)$, a ${\mathcal D}_{c}$-homothetic deformation is the following change \cite{zamkovoy}:
\begin{equation*}
\widetilde{\varphi}':=\widetilde{\varphi}, \quad \xi':=\frac{1}{c}\xi, \quad \eta':=c \eta, \quad g':=c \widetilde g + c(c-1)\eta\otimes\eta,
\end{equation*}
for some $c\neq0$.

It is known that $\mathcal{D}_c$-homothetic deformations preserve the class of paracontact metric $(\widetilde\kappa,\widetilde\mu)$-spaces. In particular, if we deform a paracontact metric $(-1,\widetilde\mu)$-space, we obtain another paracontact metric $(-1,\mu')$-space with $ \mu'=\frac{\widetilde\mu-2+2c}{c}$.

Therefore, paracontact metric $(-1,2)$-spaces remain invariant under  $\mathcal{D}_c$-homothetic deformations. Given a paracontact metric $(-1,0)$-space, if we $\mathcal{D}_c$-homothetically deform it with $c=\frac2{2-\widetilde\mu}\neq 0$ for some $\widetilde\mu\neq2$,  we will obtain a paracontact metric  $(-1,\widetilde\mu)$-space with $\widetilde\mu\neq2$. A sort of converse is also possible: given a $(-1,\widetilde\mu)$-space  with $\widetilde\mu\neq2$, a  $\mathcal{D}_c$-homothetic deformation with $c=1-\frac{\widetilde\mu}{2}\neq 0$ will give us a paracontact metric $(-1,0)$-space.

The case $\widetilde\mu=0$, $\widetilde h\neq0$  is special because the manifold satisfies \eqref{parasasakian-r} but it is not paraSasakian. Therefore, it makes sense to concentrate on $\widetilde\mu=0$ and $\widetilde\mu=2$.


We will now see that there are other possible ranks of $\widetilde h$ apart from $n$. We mention the following result, which appeared first in \cite{mio} and later with an alternative proof in \cite{mio2}.

\begin{theorem}[\cite{mio,mio2}]
Let $M$ be a $(2n+1)$-dimensional paracontact metric $(-1,\widetilde\mu)$-space $(\widetilde{\varphi},\xi,\eta,\widetilde g)$. Then we have one of the following possibilities:
\begin{itemize}
  \item
  either $\widetilde h=0$ and $M$ is paraSasakian,

  \item
  or $\widetilde h\neq 0$ and $\rank (\widetilde{h}_p)\in \{1,\ldots,n \}$ at every $p  \in M$ where $\widetilde{h}_p \neq 0$. Moreover, there exists a basis $\{ \xi_p, X_1,Y_1,\ldots,X_n,Y_n \}$ of $T_p(M)$ such that the only non-zero values of $\widetilde{g}$ on the basis are $\widetilde{g}_p(\xi_p,\xi_p)=1$ and $\widetilde{g}_p(X_i,Y_i)=\pm 1,$ and $\widetilde{h}_p$ can be written as either
      \[
      \widetilde{h}_{p| \langle X_i,Y_i \rangle}=
      \begin{pmatrix}
      0 & 0\\
      1 & 0
      \end{pmatrix}
      \quad
      \text{ or }
      \quad
      \widetilde{h}_{p| \langle X_i,Y_i \rangle}=
      \begin{pmatrix}
      0 & 0\\
      0 & 0
      \end{pmatrix},
      \]
      where obviously there are exactly $\text{rank} (h_p)$ submatrices of the first type.
\end{itemize}

If $n=1$, such a basis $\{ \xi_p, X_1,Y_1\}$ of $T_p(M)$ also satisfies that
\[
\widetilde{\varphi}_p X_1 =\pm  X_1, \quad \widetilde{\varphi}_p Y_1 =\mp  Y_1.
\]
\end{theorem}

Examples of paracontact metric $(-1,2)$-spaces with every possible constant rank of $\widetilde h$ were also presented in \cite{mio}.

\begin{example}[$(2n+1)$-dimensional paracontact metric $(-1,2)$-space with \break $\text{rank}(\widetilde h) = m \in \{1,\ldots,n \}$, \cite{mio}]
Let $\mathfrak{g}$ be the $(2n+1)$-dimensional Lie algebra with basis $\{\xi,X_1,Y_1,\ldots,X_{n},Y_{n} \}$ such that the only non-zero Lie brackets are:
\begin{align*}
[\xi,X_i]&=Y_i, \quad i=1,\ldots,m, \\
[X_i,Y_j]&=
\left\{
\begin{array}{l}
\delta_{ij} (2\xi +\sqrt2 (1+\delta_{im})Y_m)\\
\qquad +(1-\delta_{ij})\sqrt2 (\delta_{im}Y_j+\delta_{jm} Y_i),                                        \hspace*{2cm} i,j=1,\ldots,m,\\
\delta_{ij}(2\xi+\sqrt2 Y_i),                                                                          \hfill i,j=m+1, \ldots,n, \\
\sqrt2 Y_i,                                                                                             \hfill i=1, \ldots,m, \; j=m+1, \ldots,n.
\end{array}
\right.
\end{align*}

If we denote $G$ the Lie group whose Lie algebra is $\mathfrak{g}$, we can define a left-invariant paracontact metric structure $(\widetilde{\varphi},\xi,\eta,\widetilde g)$ on $G$. Indeed, let us take the $(1,1)$-tensor $\widetilde{\varphi}$ and the $1$-form $\eta$ such that
\begin{gather*}
\widetilde{\varphi} \xi=0,  \quad \widetilde{\varphi} X_i=X_i, \quad  \widetilde{\varphi} Y_i=-Y_i,  \quad i=1, \ldots,n,
\\
\eta(\xi)=1,  \quad \eta(X_i)=\eta(Y_i)=0, \quad  i=1,\ldots,n.
\end{gather*}
We define the metric $\widetilde{g}$ as the one whose only non-vanishing components are
\[
\widetilde g(\xi,\xi)=\widetilde g(X_i,Y_i)=1, \quad i=1,\ldots,n.
\]
Long computations show that the manifold is a $(-1,2)$-space and that $\text{rank}(\widetilde h)=m$.
\end{example}

Examples of $(2n+1)$-dimensional paracontact metric $(-1,0)$-spaces with  $\rank (\widetilde h)=1$ also appeared in \cite{mio} and were the first non-paraSasakian paracontact metric $(-1,\widetilde{\mu})$-spaces with $\widetilde{\mu}\neq2$ of dimension greater than $3$ that were constructed.
Later, examples of $(2n+1)$-dimensional paracontact metric $(-1,0)$-spaces with  $\text{rank}(\widetilde h)=m\in \{ 2, \ldots, n\}$ were constructed by the author in \cite{mio2}.


Finally, the question of the existence of paracontact metric $(-1,\widetilde{\mu})$-spaces with $\widetilde h$ of non-constant rank was answered also in \cite{mio2}, where the author showed the first-known examples of paracontact metric $(-1,2)$-space and $(-1,0)$-space with $\rank(\widetilde{h}_p)=0$ or $1$ depending on the point $p\in M$. We show here one of them.

\begin{example}[$3$-dimensional paracontact metric $(-1,2)$-space with $\widetilde h$ of non-constant rank, \cite{mio2}]
Let us consider the manifold $M=\mathbb{R}^3$ with the usual cartesian coordinates $(x,y,z)$. The vector fields
\[
e_1= \frac{\partial}{\partial x}+x z \frac{\partial}{\partial y}-2y\frac{\partial}{\partial z}, \quad
 e_2=\frac{\partial}{\partial y}, \quad
 \xi=\frac{\partial}{\partial z}
\]
are linearly independent at each point of $M$. We can compute
\[
[e_1,e_2]=2 \, \xi, \quad
[e_1,\xi]=-x \, e_2, \quad
[e_2,\xi]=0.
\]

We define the semi-Riemannian metric $\widetilde g$ as the non-degenerate one whose only non-vanishing components are $\widetilde g(e_1,e_2)=\widetilde g(\xi,\xi)=1$, and the $1$-form $\eta$ as $\eta=2y dx+dz$, which satisfies $\eta(e_1)=\eta(e_2)=0$, $\eta(\xi)=1$. Let $\widetilde{\varphi}$ be the $(1,1)$-tensor field defined by $\widetilde{\varphi} e_1=e_1, \widetilde{\varphi} e_2=-e_2, \widetilde{\varphi} \xi=0$. Then $\Phi=d\eta$  and $(\widetilde{\varphi},\xi,\eta,\widetilde g)$ is a paracontact metric structure on $M$.

Moreover, $\widetilde{h}\xi=0$, $\widetilde{h} e_1=x  e_2$, $\widetilde{h} e_2=0$. Hence, $\widetilde h^2=0$ and, given $p=(x,y,z) \in \mathbb{R}^3$, $\rank(\widetilde h_p)=0$ if $x=0$ and $\rank(\widetilde h_p)=1$ if $x\neq0$. Direct computations prove that the paracontact metric manifold $M$ is also a $(-1,2)$-space.
\end{example}

\end{document}